\documentclass[12pt]{amsart}
\usepackage{color}

\newcommand{\la}{\langle}
\newcommand{\ra}{\rangle}
\newcommand{\grph}{\Pi}

\newtheorem{theorem}{Theorem}[section]
\newtheorem{lemma}[theorem]{Lemma}

\newtheorem{cor}[theorem]{Corollary}
\newtheorem{proposition}[theorem]{Proposition}
\theoremstyle{definition}
\newtheorem{definition}[theorem]{Definition}

\newtheorem{remark}[theorem]{Remark}

\newtheorem{qn}[theorem]{Question}
\newtheorem{n}[theorem]{Notation}
\newtheorem{cvn}[theorem]{Convention}

\newcommand{\cg}{{\mathcal{G}}}

\newcommand{\wt}{{\mathrm{wt}}}

\newcommand{\calm}{{\mathcal M}}

\newcommand{\bfB}{(\mathbf{B})}
\newcommand{\bfS}{(\mathbf{S})}
\newcommand{\bfI}{(\textbf{I})\,}
\newcommand{\bfII}{(\textbf{II})\,}
\newcommand{\bfIII}{(\textbf{III})\,}

\newcommand{\dist}{{\mathrm{d}}}

\newcommand{\diam}{{\mathrm{diam}}\, }

\newcommand{\SL}{\mathrm{SL}}
\def\ee{\mathcal{E}}   
\def\calc{\mathcal{C}}   
\newcommand{\calh}{\mathcal{H}}
\definecolor{darkgreen}{cmyk}{1,0,1,.2}

\definecolor{m}{rgb}{1,0.1,1}

\catcode`\@=11
\long\def\@savemarbox#1#2{\global\setbox#1\vtop{\hsize\marginparwidth
  \@parboxrestore\tiny\raggedright #2}}

\newcommand {\N}{\mathbb{N}} 
\newcommand {\Z}{\mathbb{Z}}            
\newcommand {\R}{\mathbb{R}} 

\newcommand {\me}{\medskip}

\newcommand {\iv}{^{-1}}

\newcommand{\distortion}{\mathrm{dtn}}

\newcommand{\card}{\mathrm{card}\, }
\newcommand{\cale}{{\mathcal E}}
\newcommand {\fn}{\footnote}
\begin{document}

\title[{Compression functions of uniform embeddings of groups}]{Compression functions of uniform embeddings of groups into Hilbert and Banach spaces}

\author{Goulnara Arzhantseva}
\address{Universit\'{e} de Gen\`{e}ve,
Section de Math\'{e}matiques, 2-4 rue du Li\`{e}vre, Case postale
64, 1211 Gen\`{e}ve 4, Switzerland}
\email{Goulnara.Arjantseva@math.unige.ch}
\author{Cornelia Dru\c tu}
\address{University of Oxford, Mathematical Institute,
24-29 St Giles, Oxford OX1 3LB, United Kingdom.}
\email{drutu@maths.ox.ac.uk}
\author{Mark Sapir}
\address{Department of Mathematics,
1326 Stevenson Center, Vanderbilt University, Nashville, TN 37240,
USA} \email{m.sapir@vanderbilt.edu}
\thanks{This article was written while all the three authors were visitors at
the Max Planck Institute in Bonn. We are grateful to the MPI for its
hospitality. The research of the first author was supported by the
Swiss National Science Foundation, grants PP002-68627, PP002-116899.
The research of the third author was supported in part by the NSF
grants DMS 0245600 and DMS 0455881 and by a BSF (USA-Israeli)
grant.}

\subjclass[2000]{{Primary 20F69; Secondary 20F65, 20E06, 46B07}}

\date{}


\keywords{{Uniform embeddings, expander graphs, exact groups, Hilbert space
compression, asymptotic dimension, metric distortion, graph of groups.}}

\begin{abstract}
We construct finitely generated groups with arbitrary prescribed
Hilbert space compression {$\alpha\in[0,1]$}. This answers a
question of E. Guentner and G. Niblo. For a large class of Banach
spaces $\ee$ (including all uniformly convex Banach spaces), the
$\ee$--compression of these groups coincides with their Hilbert
space compression. Moreover, the groups that we construct have
asymptotic dimension {at most $2$,} hence they are exact. In
particular, the first examples of groups that are uniformly
embeddable into a Hilbert space (moreover, of finite asymptotic
dimension and exact) with Hilbert space compression 0 are given.
{These groups are also the first examples of groups with uniformly
convex Banach space compression 0}.
\end{abstract}

\maketitle

\section{Introduction}\label{prel}

\subsection{Uniform embeddings}
The property of uniform embeddability of groups into Hil\-bert spaces (and, more generally, into Banach
spaces) became popular after Gromov \cite{Gr:asinv} suggested that this property might imply the
Novikov conjecture. Indeed, following this suggestion,  Yu \cite{Yu} and later Kasparov and  Yu
\cite{KY} proved that a finitely generated group uniformly embeddable into a Hilbert space,
respectively into a uniformly convex Banach space, satisfies the Novikov conjecture.

{This raised the question whether every finitely generated group can be embedded uniformly into a
Hilbert space, or more generally, into a uniformly convex Banach space. Gromov constructed \cite{Gr}
finitely generated {\em random} groups whose Cayley graphs (quasi)-contain some infinite families of
expanders and thus cannot be embedded uniformly into a Hilbert space (or into any $\ell^p$ with $1\leq
p<\infty$, e.g.~\cite[Ch.11.3]{roe}). The recent results of V. Lafforgue \cite{LT} yield a family of
expanders that is not uniformly embeddable into any uniformly convex Banach space. Nevertheless, one
cannot apply Gromov's argument to deduce that random groups corresponding to Lafforgue's family of
graphs do not embed uniformly into any uniformly convex Banach space. Indeed, Lafforgue's expanders are
Cayley graphs of finite quotients of a non-free group, therefore there are loops of bounded size in all
of them, and the graphs have bounded girth. On the other hand, Gromov's argument succeeds only if the
girth of a graph in the family of expanders is of the same order as the diameter of the graph.


\subsection{Compression and compression gap}


\begin{definition}[cf. \cite{GuK}] Let $(X,\dist_X)$ and $(Y,\dist_Y)$ be two metric spaces and let $\phi\colon X\to
Y$ be an $1$-Lipschitz map. The {\em compression of $\phi$} is the
supremum over all $\alpha\ge 0$ such that $$\dist_Y(\phi(u),
\phi(v))\ge \dist_X(u,v)^\alpha$$ for all $u, v$ with large enough
$\dist_X(u,v)$.

If $\cale$ is a class of metric spaces closed under rescaling of
the metric, then the $\cale$\emph{--compression of }$X$ is the
supremum over all compressions of $1$-Lipschitz maps $X\to Y$,
$Y\in \cale$. In particular, if $\cale$ is the class of Hilbert
spaces, we get the {\em Hilbert space compression} of $X$.
\end{definition}

The $\cale$--compression measures the least possible distortion of distances when one tries to draw a
copy of $X$ inside a space from $\cale$. It is a quasi-isometry invariant of $X$ {and it takes values
in the interval $[0,1]$}. Similar concepts of distortion have been extensively studied (mostly for
finite metric spaces mapped into finite dimensional Hilbert spaces) by combinatorists for many years
(see \cite{Bou,De2}, for example).

Since any finitely generated group $G$ can be endowed with a word
length metric and all such metrics are quasi-isometric, one can
speak about the \emph{$\cale$--compression of a group} $G$.
Guentner and Kaminker proved in \cite{GuK} that if the Hilbert
space compression of a finitely generated group $G$ is larger than
$\frac12$ then the reduced $C^*$-algebra of $G$ is exact (in other
words, $G$ is \emph{exact} or $G$ satisfies {Guoliang} Yu's {\em
property A}~\cite{Yu}).

One of the goals of this paper is to describe all possible values
of $\cale$--compression for finitely generated groups, when
$\cale$ is either the class of Hilbert spaces or, more generally,
the class of uniformly convex Banach spaces.

A very limited information was known about the possible values of
Hilbert space compression of finitely generated groups. For example,
word hyperbolic groups have Hilbert space compression 1
\cite{BrodskiySonkin:compression}, and so do groups acting properly
and co-compactly on a cubing \cite{CN}; co-compact lattices in
arbitrary Lie groups, and all lattices in semi-simple Lie groups
have Hilbert space and, moreover, $L^p$--compression 1~\cite{Te};
any group that is not uniformly embeddable into a Hilbert space
(such groups exist by ~\cite{Gr}) has Hilbert space compression 0,
etc. (see the surveys in \cite{AGS}, \cite{Te}). The first groups
with Hilbert space compressions strictly between 0 and 1 were found
in \cite{AGS}: R. Thompson's group $F$ has Hilbert space compression
$\frac12$, the Hilbert space compression of the wreath product
$\Z\wr \Z$ is between $\frac12$ and $\frac34$ (later it was proved
in \cite{ANP} that it is actually $\frac23$), the Hilbert space
compression of $\Z \wr (\Z \wr \Z)$ is between $0$ and $\frac12$.

The notion of {$\ee$}--compression can be generalized to the
notion of {\em $\ee$--compression gap} of a space $X$, where $\ee$
is any class of metric spaces closed under rescaling of the metric
(see Definition \ref{defcomprfct}). It measures even more
accurately than the compression the best possible (least
distorted) way of embedding $X$ in a space from $\cale$. {For
example, it is proved in \cite{AGS} that R. Thompson's group $F$
has Hilbert space compression gap $(\sqrt{x}, \sqrt{x}\log x)$.
This means there exists an 1-Lipschitz embedding of $F$ into a
Hilbert space with compression function $\sqrt{x}$, and every
1-Lipschitz embedding of $F$ into a Hilbert space has compression
function at most $\sqrt{x}\log x$. This is much more precise than
simply stating that the Hilbert space compression of $F$ is
$\frac12$. Another example: it follows from \cite{Te} that every
lattice in a semi-simple Lie group has a Hilbert space compression
gap $\left( \frac{x}{\sqrt{\log x}\log\log x}, x \right)$ and the
upper bound of the gap cannot be improved. This is a much more
precise statement than the statement that the Hilbert space
compression of the lattice is $1$.}

In this paper, we show that a large class of functions appear as
Hilbert space compression functions of graphs of bounded degree, and
as upper bounds of Hilbert space compression gaps of logarithmic
size of finitely generated groups. This class of functions is
defined as follows. We use the notation $\R_+$ for the interval
$[0,\infty).$

{\begin{definition}\label{ccal} Let $\calc$ be the collection of
continuous functions $\rho :\R_+ \to \R_+$ such that {for some
$a>0$: \begin{itemize}
    \item[(1)] $\rho$ is
increasing on $[a,\infty)$, and $\lim_{x\to \infty}
\rho(x)=\infty$;
    \item[(2)] $\rho$ is subadditive;
    \item[(3)] the function
$\tau (x)=\frac{x}{\rho(x)}$ is increasing, and the function
$\frac{\tau(x)}{\log x}$ is non-decreasing on $[a, \infty)$.
\end{itemize}}
\end{definition}}

\begin{remark}\label{fctinC}
The collection $\calc$ contains all functions $x^\alpha$, for $\alpha\in (0,1)$, as well as functions
$\log x$, $\log\log x$, $\frac{x}{\log^\beta (x+1)}$, $\frac{x^\alpha}{\log^\gamma (x+1)}$ for
$\alpha\in (0,1), \beta>1, \gamma>0$, etc.
\end{remark}


\subsection{Results of the paper}

Let $\ee$ be the class of all uniformly convex Banach spaces.

\begin{proposition}[see Proposition \ref{p1}]\label{p11}
Let $\rho$ be a function in $\calc$. There exists a graph $\grph$ of bounded degree such that $\rho$ is
the Hilbert space compression function of $\grph$, and also the $\ee$--compression function of $\grph$.

In particular, for any $\alpha \in [0,1]$ there exists a graph of
bounded degree whose Hilbert space compression equals the
$\ee$--compression, and both are equal to~$\alpha$.
\end{proposition}

Using the construction of the graph in Proposition \ref{p11}, we
realize every function of $\calc$ as { the upper bound of a Hilbert
space compression gap of logarithmic size of a finitely generated
group.}

\begin{theorem}[see Theorem \ref{main}]\label{main1}
For every function $\rho \in \calc$ there exists a finitely
generated group of asymptotic dimension at most $2$ such that {for
every $\epsilon>0$, $\left(\frac{\rho}{\log^{1+\epsilon}(x+1)}\, ,\,
\rho\right)$ is a Hilbert space compression gap and
an $\ee$--compression gap of the group}.

In particular, for every $\alpha \in [0,1]$ there exists a finitely
generated group $G_\alpha$ of asymptotic dimension {at most $2$} and
with the Hilbert space compression equal to the $\ee$--compression
and equal to $\alpha$.
\end{theorem}

Since the groups $G_\alpha$ have finite asymptotic dimension, they
are all exact and uniformly embeddable\footnote{A finitely
generated group $G$ of finite asymptotic dimension has Guoliang
Yu's property A~\cite[Lemma 4.2]{HR}. That property is equivalent
to the exactness of the reduced $C^*$-algebra of $G$~\cite{Oz,HR}
and guarantees uniform embeddability into a Hilbert
space~\cite[Th.2.2]{Yu}.} into Hilbert spaces even when
$\alpha=0$. Thus we construct the first examples of groups
uniformly embeddable into Hilbert spaces, moreover exact and even
of finite asymptotic dimension, that have Hilbert space
compression $0$, and even \{uniformly convex Banach
space\}--compression zero. Note that since the construction in
\cite{Gr} does not immediately extend to uniformly convex Banach
spaces, our groups seem to be the only existing examples of groups
with \{uniformly convex Banach space\}--compression 0.

\subsection{The plan of the proofs}

The plan for proving Proposition \ref{p11} and Theorem \ref{main1} is the following. We use a family of
V.~Lafforgue's expanders $\Pi_k$, $k\geq 1$, which are Cayley graphs of finite factor-groups $M_k$ of a lattice
$\Gamma$ of $\SL_3(F)$ for a local field $F$. Lafforgue proved \cite{LT} that this family of expanders
does not embed uniformly into a uniformly convex Banach space. Now taking any function  $\rho$ in
$\calc$, we choose appropriate scaling constants $\lambda_k, k\geq 1$, such that the family of rescaled
metric spaces $(\lambda_k\Pi_k)_{k\geq 1}$ has \{uniformly convex Banach space\}-compression function $\rho$ and
Hilbert space compression function $\rho$ as well. This gives Proposition \ref{p11}.

The group satisfying the conditions of Theorem \ref{main1} is
constructed as a graph of groups. We use the fact that each $M_k$ is
generated by finitely many involutions, say $m$ (that can be
achieved by choosing a lattice $\Gamma$ generated by involutions).
One of the vertex groups of the graph of groups is the free product
$F$ of the groups $M_k$, other vertex groups are $m$ copies of the
free product $H=\Z/2\Z
* \Z$. Edges connect $F$ with each of the $m$ copies of $H$. The edge groups are free products of
countably many copies of $\Z/2\Z$. We identify such a subgroup in
$H$ with a subgroup of $F$ generated by involutions, one involution
from the generating set of each factor $M_k$. As a result, the group
$G$ is finitely generated, and each finite subgroup $M_k$ in $G$ is
distorted by a scaling constant close to $\lambda_k$. Hence the
Cayley graph of $G$ contains a quasi-isometric copy of the family of
metric spaces $(\lambda_k\Pi_k)_{k\geq 1}$. This allows us to apply
Proposition \ref{p11} and {get an upper bound for a compression gap.
A lower bound} is achieved by a careful analysis of the word metric
on $G$.

{In order to show that  $G$ has asymptotic dimension at most 2, we use a result by Dranishnikov and
Smith \cite{DS} on the asymptotic dimension of countable groups, and results by Bell and Dranishnikov
\cite{BD}, as well as by Bell, Dranishnikov and Keesling \cite{BD1} on the asymptotic dimension of
groups acting on trees.}


\subsection{Other Banach spaces}

{The class of Banach spaces to which our arguments apply cannot be extended much beyond the class of
uniformly convex Banach spaces; for instance it cannot be extended to reflexive strictly convex Banach
spaces. Indeed, our proof is based on the fact that a family of V.~Lafforgue's expanders \cite{LT} does
not embed uniformly into a uniformly convex Banach space. But by a result of Brown and Guentner
\cite{BG}, any countable graph of bounded degree can be uniformly embedded into a Hilbertian sum
$\oplus l^{p_n} (\N )$ for some sequence of numbers $p_n \in (1,+\infty )$, $p_n \to \infty $. The
Banach space $\oplus l^{p_n} (\N )$ is reflexive and strictly convex, but it is not uniformly convex.

\me

{\bf Acknowledgement.} The authors are grateful to
A.~Dranishnikov, A.~Lubotzky, A.~Rapinchuk, B.~Remy, R.~Tessera
and D.~Witte Morris for useful conversations and remarks.

\section{Embeddings of metric spaces}\label{emb}

Given two metric spaces $(X,\dist_X)$ and $(Y,\dist_Y)$ and an
$1$-Lipschitz map $\phi\colon X\to Y$ we define the
\emph{distortion} of $\phi$ \cite{LinialWigderson:expanders} as
follows:
\begin{equation}\label{eqad}
\distortion(\phi) = \max_{x\neq y} \frac{\dist_X
(x,y)}{\dist_Y(\phi(x),\phi(y))}.
\end{equation}

For a metric space $(X,\dist)$ and a collection of metric spaces
$\cale$ we define the $\cale$\emph{--distortion of} $(X,\dist)$,
which we denote by $\distortion_\cale(X,\dist )$, as the infimum
over the distortions of all $1$-Lipschitz maps from $X$ to a
metric space from $\cale$. Note that given $\lambda$ a positive
real number $\distortion_\cale(X,\dist )=
\distortion_\cale(X,\lambda \dist)$ provided that the class
$\cale$ is closed under rescaling of the metrics by $\lambda$.

\begin{remark}\label{nvertices}
If $\cale$ contains a space with $n$ points at pairwise distance
at least $1$ from each other then for every graph $X$ with $n$
vertices and edge-length metric, $\distortion_\cale X\leq \diam
X$.
\end{remark}


The notion of distortion originated in combinatorics is related to
the following notion of uniform embedding with origin in
functional analysis.

\begin{definition}\label{dunifemb}
Given two metric spaces $(X,\dist_X)$ and $(Y,\dist_Y)$, and two
proper non-decreasing functions $\rho_\pm:\R_+ \to \R_+$, {with
$\lim_{x\to \infty}\rho_{\pm} (x)=\infty$}, a map $\phi\colon X\to
Y$ is called a $(\rho_-,\rho_+)$--\emph{embedding} (also called a
\emph{ uniform embedding} or a \emph{coarse embedding}) if
\begin{equation}\label{unife}
\rho_- (\dist_X (x_1,x_2))\leq \dist_Y (\phi(x_1),\phi(x_2)) \leq
\rho_+ (\dist_X (x_1,x_2))\, {,}
\end{equation}
for all $x_1, x_2$ in $X.$


If $\rho_+(x)=Cx$, i.e. if $\phi$ is $C$-Lipschitz for some
constant $C>0$, then the embedding is called a
$\rho_-$--\emph{embedding}.
\end{definition}

\begin{definition}\label{dunifembs} For a family of metric spaces
$X_i, i\in I$, by a $(\rho_-, \rho_+)$--\emph{embedding} (resp. $\rho_-$--\emph{embedding}) of the
family we shall mean the $(\rho_-, \rho_+)$--embedding (resp. $\rho_-$--embedding) of the wedge union
of $X_i$.
\end{definition}

Let $(X,\dist_X)$ be a quasi-geodesic metric space (e.g., the set of
vertices of a graph). Then it is easy to see that any
$(\rho_-,\rho_+)$--embedding of $X$ is also a $\rho_-$--embedding.
The same holds for a family of metric spaces.


\begin{cvn}
{Since} in this paper we discuss mainly embeddings of graphs, in
what follows we restrict ourselves to {$\rho_-$--embeddings}, and
denote the function $\rho_-$ simply by~$\rho$.
\end{cvn}


\begin{n}\label{nprec}

For two functions $f,g :\R \to \R$ we write $f\ll g$ if there exist $a,b,c>0$ such that  $f(x)\leq
ag(bx)+ c $ for every $x\in \R$. If $f\ll g$ and $g\ll f$ then we write $f\asymp g$.

\end{n}

\begin{definition}\label{defcomprfct}
Let $(X,\dist )$ be a metric space, and let $\cale$ be a collection of metric spaces. Let $f,
g\colon\R_+ \to \R_+$ be two increasing functions such that $f\ll g$ and $\lim_{x\to \infty } f(x)=
\lim_{x\to\infty} g(x)=\infty .$

We say that $(f,g)$ is an  $\cale${\em--compression gap of $(X,\dist
)$ } if

\begin{itemize}
\item[(1)]
 there exists an $f$--embedding of $X$ into
a space from $\cale$;

\item[(2)]
for every $\rho$--embedding of $X$ into a space from $\cale$, we have $\rho\ll g$.

\end{itemize}

\medskip

If $f=g$ then we say that $f$ is the $\cale$--\emph{compression function} of $X$.

\medskip

The quotient $\frac{g}{f}$ is called the {\em size of the gap}. {The
functions $f$ and $g$ are called, the \emph{lower} and \emph{upper
bound of the gap} respectively.}

\medskip

The supremum of all the non-negative numbers $\alpha$
such that there is an $x^\alpha$--embedding of $X$ into a space from $\cale$ is said to be  the $\ee$--\emph{compression}
of $X$.
\end{definition}

Observe that if $\cale$ is closed under rescaling of the metrics, then
any compression gap of a metric space $X$ is a quasi-isometry
invariant.

\section{Embeddings of expanders}\label{expand}

For every finite $m$-regular graph, the largest eigenvalue of its
incidence matrix is $m$. We denote by $\lambda_2$ the second
largest eigenvalue.

We start by a well known metric property of expanders.

\begin{lemma}[\cite{lub}, Ch.1]\label{expdist}
 Let $\epsilon>0$ and let $\mathcal{G}_{m,\epsilon }$ be the
family of all $m$-regular graphs with $\lambda_2 \leq m-\epsilon$. Then there exist two constants
$\kappa ,\kappa'$ such that for any graph $\grph$ in $\mathcal{G}_{m,\epsilon}$ with set of vertices
$V$, and any vertex $x\in V$, the set $\left\{ y\in V \mid \dist(x,y)\geq \kappa \, \diam \grph
\right\}$ has cardinality at least $\kappa'\, |V|$.
\end{lemma}

We now recall properties of expanders related to embeddings into
Hilbert spaces. We denote by $\calh$ the class of all separable
Hilbert spaces.

\begin{theorem}[\cite{LLR}, Theorem 3.2(6), \cite{LinialWigderson:expanders},
Theorem 13.8 and its proof] \label{linelw}

Let $\epsilon>0$ and let $\mathcal{G}_{m,\epsilon }$ be the family of all $m$-regular graphs with
$\lambda_2 \leq m-\epsilon$.
\begin{itemize}
    \item[(i)] There exist constants $c>0$ and $d>1$ such that for any graph $\grph$ in $\mathcal{G}_{m,\epsilon}$
with set of vertices~$V$
\begin{equation}\label{distexp}
c\log |V|\le \distortion_\calh \grph \le \diam \grph \leq d\log
|V|\, .
\end{equation}

\medskip

    \item[(ii)] Moreover, there exists $r>0$ such that for every $1$-Lipschitz embedding $\phi$ of $\grph$ into a
Hilbert space $Y$ there exist two vertices $v_1$ and $v_2$ in $V$
with $\dist(v_1,v_2)\ge \kappa\, \diam \grph $ (where $\kappa$ is
the constant in Lemma \ref{expdist}) and
$$\|\phi(v_1)-
\phi(v_2)\|\le r\, .$$
\end{itemize}
\end{theorem}

\begin{remark}\label{nvertices:exp}
Thus for all expanders in $\mathcal{G}_{m,\epsilon}$ the canonical
embedding $\grph \hookrightarrow \frac{1}{\sqrt{2}}\ell^2(V)$ has
optimal distortion, see Remark~\ref{nvertices}.
\end{remark}

{V.~Lafforgue \cite{LT} proved that for some smaller family of
expanders one can replace $\calh$ by a large class of Banach
spaces.
Namely, for every prime number $p$ and  natural number $r>0$, he
defined  a class of Banach spaces
$\ee^{p,r}=\bigcup_{\alpha>0}\ee^{p,r,\alpha}$ satisfying the
following:}
\medskip

\begin{itemize}
\item[(1)] For any $p, r$, the class $\ee^{p,r}$ contains all
uniformly convex Banach spaces\fn{{A Banach space is
\emph{uniformly convex} if for every $R>0$ and every $\delta
>0$ there exists $\varepsilon = \varepsilon (R,\delta )>0$ such that
if $x,y$ are two points in the ball around the origin of radius
$R$ at distance at least $\delta $ then their middle point
$\frac{1}{2}(x+y)$ is in the ball around the origin of radius
$R-\varepsilon$.}} (including all Hilbert spaces, and even all
spaces $l^q$, $q>1$);

\medskip

\item[(2)] $\bigcup_{r} \ee^{2,r}$ is the
set of all $B$-convex Banach spaces.
\end{itemize}

\medskip

Lafforgue's family of expanders that does not embed uniformly into any Banach space from $\ee^{p,r}$ is
constructed as follows. Given numbers $p$ and $r$ as above, let $F$ be a local field such that the
cardinality of its residual field is $p^r$. Let $\Gamma$ be a lattice in $\SL(3,F)$.

The group $\Gamma$ is residually finite with Kazhdan's property
(T). {In fact,  $\Gamma$ satisfies the {\it Banach version} of property (T) with respect to the family $\ee^{p,r}$ \cite[$\S 3$ and
Proposition 4.5]{LT}.

Let $(\Gamma_k)_{k\ge 1}$ be a decreasing sequence of finite index
normal subgroups of $\Gamma$ such that $\bigcap_{k\ge 1} \Gamma_k
=\{ 1\}$, and let $M_k= \Gamma /\Gamma_k$, $k\geq 1$, be the
sequence of quotient groups.  Given a finite symmetric set of
generators of $\Gamma$ of cardinality $m\geq 2$, we consider each
$M_k$ endowed with the induced set of $m$ generators; we denote by
$\dist_k$ the corresponding word metric on $M_k$ and by $\Pi_k$
the corresponding $m$-regular Cayley graph. Since $\Gamma$ has
property (T), the family $(\Pi_k, \dist_k)_{k\geq 1}$ is a family
of expanders~\cite[Proposition 3.3.1]{lub}. Moreover, the
Banach version of property (T) for $\Gamma$ yields Proposition
\ref{propL} below, which implies that the family $(\Pi_k,
\dist_k)_{k\geq 1}$ cannot be uniformly embedded into a Banach
space from $\ee^{p,r}$ (see Remark \ref{notunife}).


\begin{proposition}[\cite{LT}, Proposition
5.2]\label{propL} For every $\alpha > 0$  there exists a constant
$C =C (\alpha)$ such that for any $k\ge 1$ and any space $Y$ from
$\cale^{p,r,\alpha}$ an $1$-Lipschitz map $\phi:\Pi_k \to Y$
satisfies:
\begin{equation}\label{ilaff}
\frac{1}{(\card \Pi_k)^2}\sum_{x,y\in \Pi_k} \|
\phi(x)-\phi(y)\|^2\leq C \, \frac{1}{\card \Pi_k} \sum_{x, y
\mbox{ \tiny{neighbors}}} \| \phi(x)-\phi(y)\|^2\, .
\end{equation}
\end{proposition}


\begin{cor}\label{Bdist}\strut

\begin{itemize}
    \item[\rm{(i)}] For every $\alpha > 0$ there exists a constant $D=D(\alpha )$ such that for any $k\geq
    1$,
\begin{equation}\label{epa}
D\, \diam \Pi_k \leq \distortion_{\cale^{p,r,\alpha}}\, \Pi_k\leq
\diam \Pi_k \, .
\end{equation}
    \item[\rm{(ii)}] Moreover, there exists $R=R(\alpha)$ such that for every $1$-Lipschitz embedding $\phi$ of $\Pi_k$
into a space $Y\in \cale^{p,r,\alpha}$ there exist two vertices
$v_1$ and $v_2$ in $\Pi_k$ with $\dist(v_1,v_2)\ge \kappa\, \diam
\Pi_k$, where $\kappa$ is the constant from Lemma \ref{expdist},
and

\begin{equation}\label{vert}
\|\phi(v_1)- \phi(v_2)\|\le R\, .
\end{equation}
\end{itemize}
\end{cor}

\proof The second inequality in (\ref{epa}) follows from Remark
\ref{nvertices} and from the fact that no Banach space can be
covered by finitely many balls of radius $1$. We prove the first
inequality in (\ref{epa}). Let $Y$ be an arbitrary space from
$\cale^{p,r,\alpha}$ and let $\phi:\Pi_k \to Y$ be an
$1$-Lipschitz map. Inequality (\ref{ilaff}) implies that
\begin{equation}\label{ineq}
\frac{1}{\left[ \distortion (\phi)\, \card \Pi_k
\right]^2}\sum_{x,y\in \Pi_k} \dist_k(x,y)^2 \leq Cm\, .
\end{equation}

Recall that each graph $\Pi_k$ is $m$-regular. Lemma \ref{expdist}
and (\ref{ineq}) imply that $Cm\, [\distortion (\phi)]^2 \geq
\kappa^2 \kappa' \left(\diam \Pi_k \right)^2$.

Now take $R>\frac{1}{\kappa} \sqrt{\frac{Cm}{\kappa'}}$, where $\kappa ,
\kappa'$ are the constants from Lemma \ref{expdist} and $C$ is the
constant from (\ref{ilaff}). Assume that there exists an
$1$-Lipschitz embedding $\phi$ of $\Pi_k$ into a space $Y\in
\cale^{p,r,\alpha}$ such that for any two vertices $v_1$ and $v_2$
in $\Pi_k$ with $\dist(v_1,v_2)\ge \kappa\, \diam \Pi_k$ the
following inequality holds:
$$
\frac{\dist(v_1,v_2)}{\|\phi(v_1)- \phi(v_2)\|}< \frac{1}{R}\, \diam \Pi_k
\, .
$$

Inequality (\ref{ilaff}) implies that
$$
Cm > \frac{1}{(\card \Pi_k)^2} \sum_{\dist (v_1,v_2) \geq \kappa
\diam \Pi_k } \frac{R^2 \dist (v_1,v_2)^2}{[\diam \Pi_k ]^2}\ge
R^2\kappa^2 \kappa' \, .
$$

This contradicts the choice of $R$.

Therefore there exist two vertices $v_1$ and $v_2$
in $\Pi_k$ with $\dist(v_1,v_2)\ge \kappa\, \diam \Pi_k$ and such that:
$$
\frac{\dist(v_1,v_2)}{\|\phi(v_1)- \phi(v_2)\|}\ge \frac{1}{R}\, \diam
\Pi_k
\, .
$$

The last inequality implies that $\|\phi(v_1)- \phi(v_2)\|\le R
\, .$
\endproof

Since $(\Pi_k, \dist_k)_{k\geq 1}$ is a family of expanders, it is
contained in some family $\cg_{m, \epsilon}$. Then according to
Theorem \ref{linelw}, we have $\diam \Pi_k \asymp \log |V_k|$.

\begin{remark}\label{notunife}
Theorem \ref{linelw}(ii) implies that no sub-family of the family
$\mathcal{G}_{m,\epsilon }$ can be uniformly embedded into a
Hilbert space, in the sense of Definition \ref{dunifembs}.

Similarly, Corollary \ref{Bdist}(ii) implies that no sub-family of
the family $(\Pi_k)_{k\ge 1}$ can be uniformly embedded into a
space from $\cale^{p,r}$.
\end{remark}

The construction in the following lemma was provided to us by Dave
Witte-Morris.

\begin{lemma}\label{witte} Let $F$ be a nonarchimedean local field of characteristic
0. There exists a lattice in the group $\SL(3,F)$ containing infinitely many noncentral involutions.
\end{lemma}

\proof Choose an algebraic number field $K$, such that

\begin{itemize}
\item $F$ is one of the (nonarchimedean) completions of $K$, and
\item $F$ is totally real.
\end{itemize}
Let $p$ be the characteristic of the residue field of $F$, and let
$c = 1 - p^3$. Thus $c < 0$, but $c$ is a square in $F$. Consider
$L= K[\sqrt{c}]$, and denote by $\tau$ the Galois involution of
$L$ over $K$. We likewise denote by $\tau$ the involution defined
on the $3\times 3$ matrices by applying $\tau$ to each entry. Let
$G = SU(3 ; L, \tau)= \left\{ g\in SL(3,L)\; ;\; \tau \left( g^T
\right)\cdot g= \mathrm{Id}_3 \right\}$. Then
\begin{itemize}
        \item[(a)] $G$ is a $K$-form of $\SL(3)$;
        \item[(b)] $G_F = \SL(3,F)$;
        \item[(c)] $G$ is compact at each real place.
\end{itemize}

Properties (a), (b), (c) can be proved using arguments similar to the ones in \cite[Chapter 10]{Wi}.

Let $S$ be the collection of all the Archimedean places of $K$ and the place corresponding to $F$.
According to the theorem due to Borel, Harish-Chandra, Behr and Harder \cite[$\S $I.3.2]{Ma}, the
$S$-integer points of $G$ form a lattice $\Gamma $ in $\SL(3,F)$. By \cite{Ta}, any lattice in a
$p$-adic algebraic group is co-compact, in particular it is the case for $\Gamma$. The lattice $\Gamma$
obviously contains noncentral diagonal matrices that are involutions.

Let $\sigma$ be one of these involutions and let $Z(\sigma )$ be
its centralizer in $\SL (3,F)$. The subgroup $\Gamma \cap Z(\sigma
)$ has infinite index in $\Gamma $. Otherwise, for some positive
integer $n$ one would have that $\gamma^n\in Z(\sigma )$ for any
$\gamma \in \Gamma$. Since $Z(\sigma )$ is an algebraic subgroup
in $\SL (3,F)$ and since $\Gamma$ is Zariski dense in $\SL (3,F)$,
it would follow that $g^n\in Z(\sigma )$ for any $g \in \SL (3,F
)$. This is impossible.

For a sequence $(\gamma_n)_{n\in \N}$ of representatives of
distinct left cosets in $\Gamma / (\Gamma \cap Z(\sigma ))$ the
involutions in the sequence $(\gamma_n \sigma \gamma_n\iv )_{n\in
\N}$ are pairwise distinct.
\endproof

\begin{lemma}\label{lmcor1}
For every prime number $p$ and  natural number $r>0$, there exists $m=m(p,r)>0$ and a lattice $\Gamma$
in $\SL(3,F)$ such that $F$ is a local field with residue field of cardinality $p^r$, and $\Gamma$ is
generated by finitely many involutions $\sigma_1,...,\sigma_m$.
\end{lemma}

\proof Let $\Gamma$ be a lattice of $\SL(3,F)$ containing infinitely many non-central involutions.
According to Lemma \ref{witte} such a lattice exists. Consider the (infinite) subgroup $\Gamma'$ of
$\Gamma$ generated by involutions. Since $\Gamma'$ is a normal subgroup in $\Gamma$, by Margulis'
Theorem $\Gamma'$ has finite index in $\Gamma$, thus it is a lattice itself.
\endproof

\begin{n}\label{nmk} Let $\Gamma$ be one of the lattices from
Lemma \ref{lmcor1}. We keep the notation given before
Proposition~\ref{propL}, referring, in addition, to the generating
set consisting of involutions. More precisely, consider a maximal
ideal $I$ in the ring of $S$-integers of the global field defining
$\Gamma$ ($S$ is the corresponding set of valuations containing
all Archimedian ones) and the congruence subgroup $\Gamma_k$ of
$\Gamma$ corresponding to the ideal $I^k$. Let
$M_k=\Gamma/\Gamma_k$, $k\ge 1$.  Let $U_k=\{ \sigma_1(k),\dots ,
\sigma_m(k) \}$ be the image of the generating set of $\Gamma$ in
$M_k$ (each $\sigma_i(k)$ is an involution). We shall denote the
corresponding word metric on $M_k$ again by $\dist_k$. Since each
$\sigma_i(k)$ is an involution, the Cayley graph $\Pi_k$ of $M_k$
is $m$-regular. Let $v(k)$ denote the cardinality of the group
$M_k$.
\end{n}

\begin{remark}\label{dm1}
{It is easy to see that $v(k)=|M_k|$ satisfies
$$ck-c_1<\log v(k)< ck+c_1$$ for some constants $c, c_1$.}
\end{remark}

\section{Metric spaces with arbitrary compression functions}

Let $(X ,\dist )$ be a metric space and $\lambda>0$. We denote by
$\lambda X$ the metric space $(X, \lambda \dist )$.

{Let $(\grph_n)_{n\geq 1}$ be the sequence of Cayley graphs of
finite factor groups $M_n$, $n\geq 1$, of the lattice $\Gamma$
from Lemma \ref{lmcor1}, see Notation \ref{nmk}.

Let $\rho$ be a function in $\calc$ (see Definition \ref{ccal}).
For every $n\geq 1$, let $x_n$ be a fixed vertex in $\grph_n$. We
are going to choose appropriate scaling constants $\lambda_n$,
$n\geq 1$, so that the wedge union $\Pi$ of the metric spaces
$\lambda_n\grph_n$, obtained by identifying all the vertices $x_n$
to the same point $x$, has the required compression function
$\rho$.



{\begin{n}\label{yn} Throughout this section we denote $\log v(n)$
by $y_n$. According to Remark \ref{dm1}, $cn-c_1<y_n< cn+c_1$ for
some constants $c, c_1$.
\end{n}}

We are looking for a sequence of rescaling constants $\lambda_n$
satisfying $\lambda_n =\rho (\lambda_n y_n)$. This is equivalent
to the fact that $\tau (\lambda_n y_n)=y_n$, where
$\tau(x)=\frac{x}{\rho(x)}$. Since $\tau(x)$ is increasing by the
definition of $\calc$, continuous and
$\lim_{x\to\infty}\tau(x)=+\infty$, $\tau\iv(x)$ exists for large
enough $x$. Thus, we can take $\lambda_n := \frac{\tau\iv
(y_n)}{y_n}$.}


\begin{proposition} \label{p1} {Let $\grph$ be the wedge union of $\lambda_n\grph_n$, $n\ge 1$.
Then $\rho$ is both the Hilbert space compression function and the
$\ee^{p,r}$--compression function of $\grph$ (up to the equivalence
relation $\asymp$).}

{In particular, $\rho $ is the \{uniformly convex Banach
space\}--compression function of $\Pi$.}
\end{proposition}

The proof of Proposition \ref{p1} will show that the same
conclusions hold if we replace the sequence ${(\lambda_n)_{n\geq
1}}$ by any sequence of numbers ${(\mu_n)_{n\geq 1}}$ satisfying
{$\mu_n \asymp \lambda_n$}.

\proof[Proof of Proposition \ref{p1}] We are going to prove that
$\rho$ is the Hilbert space compression function of $\grph$. The
proof that $\rho$ is the $\ee^{p,r}$--compression function is
essentially the same (with reference to Theorem \ref{linelw}
replaced by a reference to Corollary \ref{Bdist}).

Let $\dist_\tau$ be the canonical distance on $\grph$, and let $V$
be the set of vertices of $\grph$. Consider the map $f$ from $V$
to $\ell_2(V)$ defined as follows. For every $n\geq 1$ and every
$v\in V_n\setminus \{ x_n\}$ let $f(v)$ be equal to $\lambda_n
\delta_v$, where $\delta_v$ is the Dirac function at $v$. Also let
$f(x)=0$.

We prove that the following inequalities are satisfied, with
$d\geq 1$ the constant in Theorem~\ref{linelw}(i).

\begin{equation}\label{eq4}
\frac{1}{\sqrt{2}}\, \rho\left( \frac{\dist_\tau(v,v')}{d}\right)
\le \|f(v)-f(v')\|\le \sqrt{2}\dist_\tau(v,v')\, .
\end{equation}

Assume that $v,v'\in \grph_n$ for some $n$. Then
$\|f(v)-f(v')\|\leq \lambda_n\sqrt{2}$ (when one of the vertices
is $x_n$ the first term is $\lambda_n$). Since
$\dist_\tau(v,v')\ge \lambda_n$, we get the second inequality in
(\ref{eq4}).

To prove the first inequality, recall that by Theorem \ref{linelw}
the diameter of $\lambda_n \grph_n$ is at most $d\lambda_n\, y_n$.
Hence

$$\rho\left( \frac{\dist_\tau(v,v')}{d} \right)\le \rho \left(\lambda_n y_n\right)=\lambda_n\leq  \|f(v)-f(v')\| \, .$$

Assume now that $v\in \grph_m \setminus \{ x_m \}$ and $v'\in
\grph_n \setminus \{ x_n\}$. Then
$$
\|f(v)-f(v')\|=\sqrt{\lambda_m^2+\lambda_n^2}\leq
\lambda_m+\lambda_n\leq \dist_\tau (v,x)+\dist_\tau
(x,v')=\dist_\tau (v,v')\, .
$$

On the other hand
$$
\|f(v)-f(v')\|\geq \frac{1}{\sqrt{2}}(\lambda_m+\lambda_n)=
\frac{1}{\sqrt{2}}\left( \|f(v)-f(x)\|+\|f(v')-f(x)\| \right)\, .
$$

By the previous case, the last term is at least

$$
\frac{1}{\sqrt{2}} \left[ \rho\left( \frac{\dist_\tau(v,x)}{d}
\right) +\rho\left( \frac{\dist_\tau(x,v')}{d} \right) \right]\geq
\frac{1}{\sqrt{2}}\rho\left( \frac{\dist_\tau(v,v')}{d} \right)\,
.
$$

The latter inequality is due to the sub-additivity and the
monotonicity of $\rho$.

This proves that $\rho$ has property (1) of the Hilbert space
compression function of $\grph$, see
Definition~\ref{defcomprfct}(1).

\medskip

To prove property (2) of the Hilbert space compression function, let
us consider an increasing function $\bar{\rho}\colon \R_+\to\R_+$
with $\lim_{x\to \infty} \bar{\rho}(x)=\infty$, and an
$1$-Lipschitz $\bar{\rho}$--embedding $g$ of $\grph$ into a
Hilbert space. For any pair of vertices $v,v'$ in $\grph$ we have
\begin{equation}\label{eq13}
\bar{\rho}(\dist_\tau(v,v')) \leq \| g(v) - g(v')\|\le
\dist_\tau(v,v')\, .
\end{equation}

 We shall prove that $\bar{\rho} \ll \rho$.

For every $n\geq 1$, for any two vertices $v,v'$ in the same graph
$\lambda_n \grph_n$, the inequality (\ref{eq13}) can be re-written
as:
\begin{equation}\label{seqineq}
\bar{\rho}( \lambda_n \dist (v,v')) \leq \| g(v) - g(v')\| \leq
\lambda_n \dist (v,v')\, .
\end{equation}

We denote by $h$ the restriction of $g$ to $\lambda_n \grph_n$,
rescaled by the factor $\frac{1}{\lambda_n}$. The sequence of
inequalities (\ref{seqineq}) divided by $\lambda_n$ yields
\begin{equation}\label{corn1} \frac{1}{\lambda_n}\bar{\rho} (\lambda_n \dist(v,v')) \leq \|
h(v) -h(v')\| \leq \dist (v,v')\, .
\end{equation}

Theorem \ref{linelw} implies that for some constants $r, d$ and
$\kappa $ there exist vertices $v_1$ and $v_2$ such that
$$ \kappa y_n \leq\dist(v_1,v_2)\leq dy_n \; \; \mbox{ and }\; \; \|
h(v_1)- h(v_2)\|\le r \, .$$
 From
(\ref{corn1}) and the monotonicity of $\bar{\rho}$, we get

$$
\bar{\rho}(\kappa \lambda_ny_n)\le
\bar{\rho}(\lambda_n\dist(v_1,v_2))\le
\lambda_n \| h(v_1)- h(v_2)\|\le r\lambda_n =
r
 \rho (\lambda_n y_n).
$$

Denote by $z_n$ the product $\lambda_ny_n$, also equal to $\tau\iv (y_n)$,
by the definition of
$\lambda_n$. Then we get
\begin{equation}\label{fine}
\bar{\rho}(\kappa z_n)\le r
\rho(z_n).
\end{equation}

 We have that $y_{n}\leq cn+c_1\leq y_{n-1}+ C_1$ where $C_1=c+2c_1$.

 \medskip

 By property (3) of Definition \ref{ccal}, $\tau(x)=\frac{x}{\rho (x)}$ is
an increasing map defining a
bijection $[a,\infty ) \to [b,\infty )$. {Moreover, the condition
that $\frac{\tau (x)}{\log x}$ is non-decreasing easily implies
that for some $\theta >1$, $\tau (x) +1 \leq \tau (\theta x)$ for
every $x\geq a$.} It follows that  $\tau (x) +C_1 \leq \tau (C x)$
for every $x\geq a$, where $C$ is a power of $\theta$ depending on
$C_1$. This implies that for every $y\geq b$, $\tau\iv(y+C_1)\leq
C \tau\iv(y)$. Consequently, for $n$ large enough, $z_n =\tau\iv
(y_n)\leq \tau\iv (y_n+C_1)\leq C \tau\iv (y_{n-1})=z_{n-1}$. We
thus get
$$\frac{z_n}{z_{n-1}}\leq C\, .$$

Now take any sufficiently large $x>0$. Then for some $n$, $x$ is
between $z_{n-1}$ and $z_n$. Hence (by the monotonicity of
$\bar\rho$ and $\rho$), we get:

$$\bar\rho(\kappa x)<\bar\rho(\kappa z_n)\le r\rho(z_n)
<r
\rho(Cz_{n-1})
\ll \rho(Cx).$$ Therefore, $\bar\rho\ll \rho$ as required.
\endproof

One can obviously replace $\grph$ in Proposition
\ref{p1} by a uniformly proper space (graph of bounded degree)
with the same property.

The following corollary immediately follows from Proposition
\ref{p1}.

\begin{cor} For every number $\alpha$ in $[0,1]$ there
exists a proper metric space (graph of bounded degree) whose
Hilbert space compression and the $\ee^{p,r}$--compression are
equal to $\alpha$.
\end{cor}

{Note that $\alpha=0$ can be obtained by taking, say, $\rho (x)= \log x$ in Proposition \ref{p1}. To
achieve $\alpha=1$ take the one-vertex graph.}

\section{Discrete groups with arbitrary Hilbert space compressions}\label{gr}

For every prime number $p$ and every natural number $r>0$, we
denote $\ee^{p,r}$ simply by $\ee$ in what follows. Recall that
$\ee$ contains all the uniformly convex Banach spaces.

Pick a function $\rho\in\calc$. For simplicity we assume that $\rho
(1)>0$. {As in Definition \ref{ccal}, we denote by $\tau$ the
function $\frac{x}{\rho (x)}$.}

Let $(M_k, \dist_k)_{k\ge 1}$ be the sequence of finite groups
 defined as above, for fixed $p$ and $r$, see Lemma~\ref{lmcor1} and Notation~\ref{nmk}.
Recall that $\dist_k$ is the word metric associated to the
generating set $U_k=\{\sigma_1(k),\dots ,\sigma_{m}(k)\}$
consisting of $m=m(p,r)$ involutions, and that $v(k)$ is the
cardinality of $M_k$.

\begin{n}\label{notk}
We denote by $|\cdot |_k$ the length function on $M_k$ associated
to $U_k$.
\end{n}

\begin{n}\label{mmu}
Let us fix three sequences of numbers: $(\lambda_k)_{k\ge 1}$,
$(m_k)_{k\geq 1}$, and $(\mu_k)_{k\geq 1}$.

For $k\geq 1$, we set $\lambda_k$ such that $\rho (\lambda_k \log v(k))=\lambda_k$; equivalently
$\lambda_k =\frac{\tau\iv (\log v(k))}{\log v(k)}$ with $v(k)$ as above. Without loss of generality we
assume (by taking a suitable subsequence of $M_k$) that for every $k\geq 1$, $\lambda_{k+1}-\lambda_k
\geq 4$.

For every $k\ge 2$, let $m_k$ be the integer part of
$\frac{\lambda_k -1}{2}$, and let $m_1=0$.  Then we put $\mu_k
=2m_k+1$ for every $k\geq 1$.
\end{n}

According to our choice, for every $k\geq 2$
\begin{equation}\label{nk}
\mu_k \le \lambda_k < \mu_k+2\, .
\end{equation}

We are now ready to construct our group.

Let $F$ be the free product of the groups $M_k$, $k\ge 1$.

For every $i\in \{1,\dots, m\}$, consider a copy $H_i$ of the free
product $\Z/2\Z * \Z$, where the generator of the $\Z/2\Z$-factor
of $H_i$ is denoted by $\sigma_i$, while the generator of the
$\Z$-factor of $H_i$ is denoted by $t_i$.

\begin{n}\label{hi}
We denote by $|\cdot |_{H_i}$ the length function on $H_i$
relative to the generating set $\{ \sigma_i, t_i\}$. For every
$k\in\Z$, we denote the element $t^k_i\sigma_i t^{-k}_i$ in $H_i$
by $\sigma^{(k)}_i$. Note that $| \sigma^{(k)}_i |_{H_i}=2k+1$,
and that $H_i$ is the semidirect product of $\la \sigma^{(k)}_i,
k\in \Z {\ra} $ and $\la t_i \ra$.
\end{n}

Let $G$ be the fundamental group of the following graph of groups
(see \cite[Ch.5.1]{ser} for the definition). The vertex groups are
$F$ and $H_1,...,H_{m}$, {the only edges of the graph are $(F,
H_i)$, $i=1,...,m$.} The edge groups $F\cap H_i$ are free products
{\Large{$\ast$}}$_{k\geq 1}\left( \Z /2\Z \right)_{k}$ of
countably many groups of order $2$, where the $k$-th factor
$\left( \Z /2\Z \right)_{k}$ is identified with $\la
\sigma_i(k)\ra < M_k $ in $F$ and with $\la\sigma_i^{(m_k)}\ra $
in $ H_i$.

\begin{n}\label{dg}
We denote by $\dist_G$ and by $|\cdot |_G$ the word metric and
respectively the length function on $G$ associated to the
generating set $U=\{ \sigma_1(1),...,\sigma_{m}(1),
t_1,...,t_{m}\}$.
\end{n}

\begin{theorem}\label{main}
With the notation above and the terminology in Definition
\ref{defcomprfct}, the following hold:

\begin{itemize}
\item[\rm{\bfI}] {The group $G$ has as a Hilbert space compression gap, and an
$\ee$--compression gap $(\sqrt{\rho} ,\rho )$.}

\item[\rm{\bfII}] {The group $G$ has as a Hilbert space compression gap, and
an $\ee$--compression gap $\left( \frac{\rho}{\log^{1+\epsilon}(x)},\rho \right)$, for every
$\epsilon>0$.} Hence the Hilbert space compression and the $\ee$--compression of $G$ are equal to the
supremum over all the non-negative numbers $\alpha$ such that $x^\alpha \ll \rho$.

\item[\rm{\bfIII}] The asymptotic dimension of $G$ is {at most
$2$.}



\end{itemize}


\end{theorem}

We start our proof by describing the length function $|\cdot |_G$.

By the standard properties of amalgamated products \cite[Ch.IV,
Th.2.6]{LyndonSchupp}, the free product $F= \mbox{\Large $\ast$
}_{k\geq 1} M_k$ and the groups $H_i$, $i\in \{1,2,\dots ,m\}$,
are naturally embedded into $G$. Moreover, the subgroup $T=\la
t_i, 1\le i\le m\ra$ is free of rank $m$ and a retract of $G$, so
its length function coincides with $|\cdot|_G$ restricted to $T$.

\begin{definition}\label{weight} With every element $s\in M_k$ we
assign the \emph{weight} $\wt(s)=\mu_k\, |s|_k$. With every
element $b\in H_i$ we assign the \emph{weight} $\wt(b)=|b|_{H_i}$.
The\emph{ weight} $\wt(W)$ of any word $W=s_1b_1...s_nb_n$ with
$s_1$ and $b_n$ possibly trivial and $n\in \N$ is the sum
$\sum_{i=1}^n [\wt(s_i) + \wt(b_i)]$.
\end{definition}



Let $F'$ be the free product of $F=\mbox{\Large $\ast$ }_{k\geq 1}
M_k $ and all the $2$-element groups $\left\la \sigma_i^{(k)}
\right\ra $ with $k\in \Z \setminus \{m_k \; {,}\; k\geq 1\}$ and
$i\in \{1,\dots, m\}$. Then $G$ is the multiple HNN extension of
$F'$ with free letters $t_1,\ldots,t_m$ shifting $\sigma_i^{(k)}$.
This is the presentation of $G$ used in the following lemma.

A subword $t_i^{\pm 1}\sigma_i^{(k)}t_i^{\mp 1}$ of a word in
generators $\sigma_i^{(k)}, t_i$ of $G$ is called a \emph{pinch}.
It can be removed and replaced by $\sigma_i^{(k\pm 1)}$; we call
this operation \emph{removal of a pinch}.

\begin{lemma}\label{lengthofk} For every $k\geq 1$ and $g\in M_k$,
$|g|_G=\wt(g)=\mu_k \, |g|_k$.
\end{lemma}

\proof Let $s=|g|_k$ and let
$g=\sigma_{i_1}(k)\sigma_{i_2}(k)...\sigma_{i_s}(k)$ be a shortest
representation of $g$ as a product of generators of $M_k$. Then we
can represent $g$ as
$$t_{i_1}^{m_k}\sigma_{i_1}(1)t_{i_1}^{-m_k}
t_{i_2}^{m_k}\sigma_{i_2}(1)t_{i_2}^{-m_k}...t_{i_s}^{m_k}\sigma_{i_s}(1)t_{i_s}^{-m_k}
.$$ Hence $|g|_G\le (2m_k+1) \, s=\mu_k \, |g|_k$.

Let $W$ be any shortest word in the alphabet $U$ representing $g$
in $G$. Since $g\in F$, there exists a sequence of removals of
pinches and subwords of the form $aa\iv$ that transfers $W$ into a
word $W'$ in generators of $F$ (i.e. letters of the form
$\sigma_j(l)$) representing $g$. Since the weight of every
$\sigma_i(1)$ and $t_i^{\pm 1}$ is $1$, the total weight of the
word $W$ is the length of $W$. The removal of a pinch does not
increase the total weight, and the removal of a word $aa\iv$
decreases it. Hence the total weight of the letters in $W'$ is at
most $|W|$. On the other hand, since $W'$ is a word in generators
of the free product $F$ representing an element of one of the
factors $M_k$, all letters from $W'$ must belong to $M_k$, i.e.
they have the form $\sigma_i(k)$ with $i\in\{1,\dots, m\}$.
Therefore the total weight of $W'$ is $\mu_k \, |W'|$. Since $W'$
represents $g$ in $M_k$, we have that the total weight of $W'$ is
at least $\mu_k\, |g|_k$. Hence $|g|_G=|W|\ge \mu_k\, |g|_k$ as
required.
\endproof

The choice of $\mu_k$ given by (\ref{nk}),  Proposition \ref{p1}
and the construction in its proof, and Lemma \ref{lengthofk} imply
the following.

\begin{lemma}\label{comprs}\strut

\begin{itemize}
    \item[\rm{(1)}] The function $\rho$ is the $\cale$--compression function and the
Hilbert space compression function of $S=\bigcup_{k\geq 1}M_{k}$
endowed with the restriction of the metric $\dist_G$.
    \item[\rm{(2)}] The map $f: S \to l_2 (S)$ defined by $f(1)=0$ and
    $f(s)=\mu_k \delta_s$ for every $s\in M_k \setminus \{1\}$ is
    a $\rho$--embedding of $(S, \dist_G)$.
\end{itemize}
\end{lemma}

Lemma \ref{comprs} immediately implies

\begin{lemma}\label{comprs1}{The function $\rho$ is an upper bound
for an $\cale$--compression gap of $G$ (in particular, also for a
Hilbert space compression gap of $G$).}
\end{lemma}

\begin{definition}\label{nf}
We are going to use the standard normal forms of elements in the
amalgamated product $G$ (see \cite[Ch.IV.2]{LyndonSchupp}):

\begin{equation}\label{eqnf}
s_1b_1...s_nb_n
\end{equation}
where

\begin{itemize}
\item[$(N_0)$] possibly $s_1=1$ or $b_n=1$; \item[$(N_1)$] $s_i\in
M_{k_i}$ and $b_i\in H_{l_i}$ ($s_i, b_i$ are called {\em
syllables}); \item[$(N_2)$] if $b_i=1$, $i<n$, then $s_i$ and
$s_{i+1}$ are not in the same $M_k$; if $s_i=1$, $i>1$, then
$b_{i-1}, b_i$ are not in the same $H_k$.
\end{itemize}
\end{definition}

One can get from one normal form to another by sequences of
operations as follows:
\begin{enumerate}
\item[($r_1$)]replacing $s_ib_i$ with $(s_i\sigma_k(p))
(\sigma_k(p)b_i)$ ({\em left insertion});

\item[($r_2$)] replacing $b_is_{i+1}$ with
$(b_i\sigma_k(p))(\sigma_k(p)s_{i+1}{)}$ ({\em right insertion});

\item[($r_3$)] replacing two-syllable words with syllables from
the same factor by the one-syllable word, their product.
\end{enumerate}

One can easily check that this collection of moves is confluent,
so one does not need the inverse of rule ($r_3$).

Choose one representative for each left coset of $M_k/\la
\sigma_i(k)\ra$ and one representative in each left coset of
$H_i/\la \sigma_i(k)\ra$, the representative of $\la
\sigma_i(k)\ra$ is $1$. Let $\calm$ be the set of all these
{representatives}. We consider only the normal forms (\ref{eqnf})
in which all the $s_i$ and $b_i$ (except possibly for the last
non-identity syllable) are in $\calm$. We shall call these normal
forms {\em good}. Every element of $G$ has a unique good normal
form.

\begin{lemma} \label{lengthinG} Let $V=s_1b_1s_2...s_nb_n$ be the
good normal form of an element $g$ in the amalgamated product $G$.
Then the length $|g|_G$ is between $\wt(V)/3$ and $\wt(V)$, i.e.
\begin{equation}\label{eq45}
\sum_{i=1}^n \mu_{k_i}|s_i|_{{k_i}}+\sum_{i=1}^n
|b_i|_{H_{l_i}}\ge |g|_G\ge \frac13{\left(\sum_{i=1}^n \mu_{k_i}
|s_i|_{k_i}+\sum_{i=1}^n |b_i|_{H_{l_i}}\right)}.\end{equation}
For every normal form $V'$ of $g$, the length of $g$ is between
$\wt(V')/9$ and $\wt(V')$.
\end{lemma}

\proof {The first inequality in (\ref{eq45}) follows from Lemma
\ref{lengthofk} and it holds for every normal form.}

Consider a shortest word $W$ in the generators of $G$ representing
$g$. Making moves of type $(r_3)$, we rewrite $W$ in a normal form
$W'$ without increasing its weight (note that a removal of pinches
can be seen as a succession of two  $(r_3)$-type moves). We have
that $|g|_G= \wt (W)\geq \wt (W') \geq |g|_G$; the last inequality
holds because any word of the form (\ref{eqnf}) and of weight
$\ell$ is equal in $G$ to a word of length $\ell$ in the
generators in $U$ and their inverses. It follows that $\wt
(W')=|g|_G$. By doing insertion moves $(r_1)$ and $(r_2)$, we can
rewrite the normal form $W'$ into the good normal form $V$ of $g$.
Note that each syllable is multiplied by at most two involutions
during the process. The weight of each of these involutions does
not exceed the weight of the syllables $s_i, s_{i+1}$ involved in
the moves. Hence the total weight of the word cannot more than
triple during the process. Hence $\wt(V)\le 3\wt(W')=3|g|_G$ which
proves the second inequality in (\ref{eq45}).

The second statement is proved in a similar fashion: one needs to
analyze the procedure of getting the good normal form from any
normal form, and then use the first part of the lemma.
\endproof

\begin{definition}\label{ext}
Let $s_1b_1...s_nb_n$ be the good normal form of $g$. Represent
each $b_i$ as the shortest word $w(b_i)$ in the alphabet of
generators {$\left\{ \sigma_{l_i}, t_{l_i} \right\}$ of the
corresponding subgroup $H_{l_i}$}. Then the word
$s_1w(b_1)...s_nw(b_n)$ in the alphabet $\bigcup_{k\ge 1} M_k
\cup{ \{\sigma_i, t_i, i\in \{1,2,...,m\}\}}$ is called the {\em
extended normal form} of $g$.
\end{definition}

{The unicity of the good normal form and of each word $w(b_i)$
implies that every element in $G$ has a unique extended normal
form.}

\begin{definition}\label{pqf}
For every pair of elements $g$ and $h$ in $G$ with extended normal
forms $g=s_1w(b_1)...s_lw(b_l)$ and $h=s_1'w(b_1')...s_n'w(b_n')$
let $p(g,h)=p(h,g)$ be the longest common prefix of these words,
of length $i(g,h)$. Thus
$$g\equiv
p(g,h)f(g,h)q(g,h)\;\; \; \mbox{ and }\; \;  h\equiv
p(g,h)f(h,g)q(h,g)
$$ where the words
$f(g,h)q(g,h), f(h,g)q(h,g)$ have different first syllables
$f(g,h)$ and $f(h,g)$ respectively. Here the \emph{syllables} are
either in some $M_k\, ,\, k\geq 1$, or they are in $ \{t_1^{\pm
1},t_2^{\pm 1},..., t_{m}^{\pm 1}\}\, .$

Note that $p(g,h), q(g,h), q(h,g)$ are extended normal forms;
$f(g,h)$ (resp. $f(h,g)$) is either an element in $M_k$ for some
$k\geq 1$ or it is in {$\{\sigma_i, t_i^{\pm 1}\}$} for some $i\in
\{1,2,...,m\}$.
\end{definition}

\begin{definition}\label{bar} Let $s_1w(b_1)...s_nw(b_n)$ be the extended
normal form of $g$.  Then for every $i$ let $g[i]$ be the $i$-th
letter of the extended normal form, let $\hat g_i$ be the prefix
ending in $g[i]$, and let $\bar g_i$ be the suffix starting with
$g[i]$ of the extended normal form.
\end{definition}

\begin{lemma}\label{distg} For every $g, h\in G$, the distance
$\dist_G(g,h)$ is in the interval $[A/9, A]$ where $A$ is equal
to:
\begin{enumerate}
\item[$\bfS$] $\mu_k\, \dist_k(f(g,h),
f(h,g))+\wt(q(g,h))+\wt(q(h,g))$ if $f(g,h), f(h,g)\in M_k$, or

\item[$\bfB$] $\wt(f(g,h)q(g,h))+\wt(f(h,g)q(h,g))$ otherwise.
\end{enumerate}
\end{lemma}

\proof {In Case $\bfS$, $q(g,h)\iv (f(g,h)\iv f(h,g))q(h,g)$
becomes a normal form for $g\iv h$ if we combine all neighbor
letters from the same $H_k$ into one syllable, and $f(g,h)\iv
f(h,g)$ into one syllable. In Case $\bfB$, $q(g,h)\iv f(g,h)\iv
f(h,g)q(h,g)$ becomes a normal form for $g\iv h$ after a similar
procedure. Then one can use Lemma \ref{lengthinG}.\endproof}

For every element $g$ in $G$ whose extended normal form has the
last syllable in $Y$, where $Y=M_k$ for some $k\geq 1$ or
$Y=\{\sigma_i , t_i^{\pm 1}\}$ for some $i\in\{1,\dots, m\}$, we
consider a copy $\phi_g$ of an $1$-Lipschitz embedding of $Y$ into
a Hilbert space $\calh_g$ with optimal distortion: either the
embedding from Remark \ref{nvertices:exp} rescaled by the factor
$\mu_k$ if $Y=M_k$ (i.e. the embedding sending all {non-trivial}
elements in $M_k$ in pairwise orthogonal vectors of length
$\mu_k$) or the embedding defined by any choice of an orthonormal
basis in a copy of $\R^2$ otherwise. Note that since $\mu_1=1$
this is coherent with the identification of $\sigma_i\in H_i$ with
$\sigma_i(1)\in M_1$.

Let $\calh$ be the Hilbertian sum of all $\calh_g$. Consider the
following map from $G$ into $\calh$: for every $g\in G$ let

\begin{equation}\label{psi}
\psi(g):=\sum_i \phi_{\hat g_i}(g[i]).
\end{equation}

We use the map $\psi$ to prove the following.

\begin{lemma}\label{comp1}
{The group $G$ has as a Hilbert space and an $\ee$--compression
gap $(\sqrt{\rho },\rho )$.}
\end{lemma}

\proof In view of Lemma \ref{comprs1}, it suffices to show that
$\psi$ is a $\sqrt{\rho}$--embedding to finish the proof.

Let $g, h\in G$. Then, by (\ref{psi}), since any two elements of
$\calh$ of the form $\phi_{\hat g_i}(g[i])$, $\phi_{\hat
h_j}(h[j])$ either coincide or are othogonal to each other, we
have:

$$\begin{array}{llllll}
\|\psi(g)-\psi(h)\|^2& = & \|\phi_{p(g,h)f(g,h)}(f(g,h))\|^2& + &
\|\phi_{p(g,h)f(h,g)}(f(h,g))\|^2 &
\\[3pt]
& + &\sum_{j>i(g,h)+1} \|\phi_{\hat g_j}(g[j])\|^2 & +&
\sum_{j>i(g,h)+1}\|\phi_{\hat h_j}(h[j])\|^2 & \le
C\dist_G(g,h)^2\end{array}
$$
for some constant $C$ (by Lemma \ref{distg} and the fact that all
embeddings $\phi$ are $1$-Lipschitz). Thus $\psi$ is a Lipschitz
embedding.

To prove the lower bound, note that for $g[i]$ in $H_k$,
$$|g[i]|_G=\|\phi_{\hat g_i}(g[i])\|=1\gg \rho(|g[i]|_G).$$
Also if $g[i]\in M_k$ then by Lemma \ref{comprs} (2), and by the
definition of $\phi_g$ we may write
$$|g[i]|_{G}=\mu_k\, |g[i]|_{k}=\|\phi_{\hat g_i}(g[i])\|\gg \rho(|g[i]|_G)\, .$$

Thus, in all cases, $$\| \phi_{\hat g_i}(g[i])\|\gg
\rho(|g[i]|_G).$$

Now, {in} case $\bfB$,
$$\|\psi(g)-\psi(h)\|^2\gg \rho(\dist_G(f(g,h),f(h,g)))^2+
\sum_{j>
i(g,h)+1}\rho(|g[j]|_G)^2+\sum_{j>i(g,h)+1}\rho(|h[j]|_G)^2.$$ The
fact that $\rho(g[i])\ge \rho(1)$ for non-identity $g[i]$, and the
subadditivity of $\rho$ (plus Lemma \ref{distg}) imply
$$\|\psi(g)-\psi(h)\|^2\gg \rho(\dist_G(g,h))$$ as desired.
In case $\bfS$, the proof is similar, only one needs to take into
account that
$$\begin{array}{lll}\|\phi_{p(g,h)f(g,h)}(f(g,h))\|^2 & +\,\,\,
\|\phi_{p(g,h)f(h,g)}(f(h,g))\|^2 & \\
& =\| \phi_{p(g,h)f(g,h)}(f(g,h)) -
\phi_{p(g,h)f(h,g)}(f(h,g))\|^2 &
\\ & \ge \rho(\dist_G(f(g,h),f(h,g)))^2.\end{array}$$
\endproof

Lemma \ref{comp1} gives part \bfI of Theorem \ref{main}.


Let us prove part \bfII. For simplicity, we take $\epsilon=1$. The
reader can easily modify the proof to make it work for every
$\epsilon>0$. {Thus, we are going to prove that
$\frac{\rho(x)}{\log^2(x+1)}$ is a lower bound for a Hilbert space
compression gap and an $\ee$--compression gap of $G$. That is, we
are going to prove that there exists a
$\frac{\rho(x)}{\log^2(x+1)}$--uniform embedding of $G$ in a Hilbert
space.}

Consider the same Hilbert space $\calh$ and the same embeddings
$\phi_g$ for $g\in G$ as before. Let us define an embedding $\pi$
of $G$ into $\calh$. Let $g$ be given in an extended normal form.
Then we set

$$\pi(g):=\sum_j \kappa_j \phi_{\hat g_j}(g[j])$$
where the coefficients $\kappa_j$ are defined as follows:

$$\kappa_j=\frac{\sqrt{\wt(\bar g_j)}}{\log (\wt(\bar g_j)+1)\sqrt{\wt(g[j])}}.$$

We shall need the following two elementary inequalities, the first
of which is obvious.

\begin{lemma}\label{new} For every sequence of positive numbers $a_1,...,a_n$, we have
$$\sum_{i=1}^n \left(\sum_{r=i}^n a_r\right)a_i \le
\left(\sum_{r=1}^n a_i\right)^2.$$
\end{lemma}


\begin{lemma}
\label{lemineq} Let $a_1,...,a_n$ be positive real numbers with
$a_n\ge 1$. Then
\begin{equation}\label{eq34}
\begin{array}{l}
\frac{a_1}{(a_1+...+a_n)\log^2 (a_1+...+a_n+1)}+
\frac{a_2}{(a_2+...+a_n)\log^2
(a_2+...+a_n+1)}+...+\frac{a_n}{a_n\log^2 (a_n+1)}\le C.
\end{array}
\end{equation}
for some constant $C$.
\end{lemma}

\proof For every $k=1,...,n$ denote $s_k=a_k+...+a_{n}$. Then
(\ref{eq34}) can be rewritten as follows:

\begin{equation}\label{eq35}
\frac{s_1-s_2}{s_1\log^2 (s_1+1)}+...+
\frac{s_{n-1}-s_n}{s_{n-1}\log^2 (s_{n-1}+1)}+ \frac{1}{\log^2
(s_n+1)}.
\end{equation}

Since the function $\frac{1}{x\log^2 (x+1)}$ is decreasing, one
can estimate (provided $b\geq a\geq 1$)
$$\frac{b-a}{b\log^2(b+1)}\le \int_a^b
\frac{dx}{x\log^2(x+1)}\le 2\int_a^b
\frac{dx}{(x+1)\log^2(x+1)}=\frac{2}{\log (a+1)}-\frac{2}{\log
(b+1)}.$$

Applying this inequality to (\ref{eq35}) and using $a_n\geq 1$, we
get (\ref{eq34}).
\endproof

\me

Let $g, h$ be two elements in $G$ given in their extended normal
forms. Let $i=i(g,h)$, $f=f(g,h)$, $f'=f(h,g)$, $p=p(g,h)$,
$q=q(g,h)$, $q'=q(h,g)$.

Then {$\pi(g)-\pi(h)=\Sigma_1+\Sigma_2+\Sigma_3$} where

$$\Sigma_1=\sum_{r=1}^i \left( \frac{\sqrt{\wt ({\bar g}_r)}}{\log (\wt ({\bar g}_r) +1)} -
\frac{\sqrt{\wt ({\bar h}_r)}}{\log (\wt ({\bar h}_r)
+1)}\right)\frac{1}{\sqrt{\wt (g[r])}}\phi_{\hat g_r}(g[r]),$$

\bigskip

$$\Sigma_2=\frac{\sqrt{\wt(\bar g_{i+1})}}{\log (\wt(\bar g_{i+1})+1)
\sqrt{\wt(g[i+1])}} \phi_{pf}(f)-\frac{\sqrt{\wt({\bar
h_{i+1}})}}{\log (\wt(\bar h_{i+1})+1)\sqrt{\wt(h[i+1])}}
\phi_{pf'}(f'),$$

\bigskip

$$\Sigma_3=\sum_{{r}\ge i+2}
\frac{\sqrt{\wt(\bar g_r)}}{\log (\wt(\bar
g_{r})+1)\sqrt{\wt(g[r])}}\phi_{\hat g_r}(g[r])- \sum_{{r}\ge i+2}
\frac{\sqrt{\wt(\bar h_r)}}{\log (\wt(\bar
h_{r})+1)\sqrt{\wt(h[r])}}\phi_{\hat h_r}(h[r]).$$

\bigskip

\begin{lemma}[\bf Upper bound] \label{upperbound}
$\|\pi(g)-\pi(h)\|\ll \dist_G(g,h)$, that is $\pi$ is Lipschitz.
\end{lemma}

\proof Since $G$ is finitely generated, it suffices to prove that
the norm $\|\pi(g)-\pi(gs)\|$ is bounded uniformly in $s\in U \cup
U\iv$ and $g\in G$. In this case the sum $\Sigma_3$ does not
appear, and $\Sigma_2= \frac{1}{\log 2}\phi_{gs}(s)$ has norm
$\frac{1}{\log 2}$. By eventually replacing $g$ with $gs$ and $s$
with $s\iv$, we can always assume that $h=gs$ satisfies $\wt(\bar
h_r)=\wt(\bar g_r)+1$ for all $r\leq i$. Then

$$
\|\Sigma_1\|^2=
\sum_{r=1}^i\left(\frac{\sqrt{\wt(\bar{g}_r)+1}}{\log
(\wt(\bar{g}_r)+2)}-\frac{\sqrt{\wt(\bar{g}_r)}}{\log
(\wt(\bar{g}_r)+1)} \right)^2\, \wt (g[r])\, .
$$

Since the function $x\mapsto \frac{\sqrt{x}}{\log (x+1)}$ is
increasing for $x\geq e^2 -1$ it follows that

$$
\|\Sigma_1\|^2 \ll \sum_{r=1}^i\frac{\left[
\sqrt{\wt(\bar{g}_r)+1}-\sqrt{\wt(\bar{g}_r)}\right]^2}{\log^2
(\wt(\bar{g}_r)+1)} \, \wt (g[r])\leq \sum_{r=1}^i
\frac{\wt(g[r])}{\wt(\bar g_r)\log^2(\wt(\bar g_r)+1)}.
$$

Lemma \ref{lemineq} implies now that
\begin{equation}\label{sigma1}\|\Sigma_1\|^2= O(1).\end{equation}
\endproof

\begin{lemma}[\bf Lower bound]\label{lowerbound} $\|\pi(g)-\pi(h)\|\gg
\frac{\rho(\dist_G(g,h))}{\log^2(\dist_G(g,h)+1)}$ if $g\ne h$.
\end{lemma}

\proof We can write
\begin{equation}\label{lowers3}
\begin{array}{ll}
\left\| \Sigma_3 \right\|^2 & \gg \sum_{r\ge i+2} \frac{\wt(\bar
g_r) }{\log^2(\wt(\bar g_r)+1)\wt(g[r])} \rho \left(\wt(g[r])
\right)^2
\\
\\
& +\sum_{r\ge i+2} \frac{\wt(\bar h_r) }{\log^2(\wt(\bar
h_r)+1)\wt(h[r])} \rho \left(\wt(h[r]) \right)^2\, .
\end{array}
\end{equation}

If $g,h$ are in Case $\bfB$ then
$$
\left\| \Sigma_2 \right\|^2 \gg \frac{\wt(\bar g_{i+1}
)}{\log^2(\wt(\bar g_{i+1})+1)\wt(f)} \rho \left(\wt(f) \right)^2
+ \frac{\wt(\bar h_{i+1})}{\log^2(\wt(\bar h_{i+1})+1)\wt(f')}
\rho \left( \wt(f')\right)^2\, .
$$

Thus in case $\bfB$ we have
$$
\begin{array}{ll}
\left\| \Sigma_2 \right\|^2 + \left\| \Sigma_3 \right\|^2 & \gg
\sum_{r>i} \frac{\wt(\bar g_r)}{\log^2(\wt(\bar g_r)+1)\wt(g[r])}
\rho \left( \wt(g[r])\right)^2 \\
\\
& + \sum_{r>i} \frac{\wt(\bar h_r) }{\log^2(\wt(\bar
h_r)+1)\wt(h[r])} \rho \left(\wt(h[r]) \right)^2\, .
\end{array}
$$

By the subadditivity of $\rho$, we have
$$\mathbf{F}= \left[
\sum_{r>i} \rho (\wt(g[r])) + \sum_{r>i} \rho
(\wt(h[r]))\right]^2\ge \rho(\dist_G(g,h))^2.$$

 By the Cauchy-Schwartz
inequality, the {squared} sum $\mathbf{F}$ does not exceed the
product $F_1F_2$ where
$$
F_1=\sum_{r>i} \frac{\wt(g[r])\log^2(\wt(\bar g_r)+1)}{\wt(\bar
g_r)} + \sum_{r>i} \frac{\wt(h[r])\log^2(\wt(\bar
h_r)+1)}{\wt(\bar h_r)}
$$

$$F_2=\sum_{r>i} \frac{\wt(\bar g_r)}{\log^2(\wt(\bar g_r)+1)\wt(g[r])}
\rho(\wt(g[r]))^2 + \sum_{r>i} \frac{\wt(\bar
h_r)}{\log^2(\wt(\bar h_r)+1)\wt(h[r])} \rho(\wt(h[r]))^2 \, .
$$

Then $\frac{\rho(\dist_G(g,h))^2}{\log^4(\dist_G(g,h)+1)}$ does
not exceed $\frac{1}{\log^4(\dist_G(g,h)+1)}\mathbf{F}$ that, in
turn, does not exceed

$$
\frac{1}{\log^4(\dist_G(g,h)+1)} F_1 F_2\leq \left( \sum_{r>i}
\frac{\wt(g[r])}{\log^2(\wt(\bar g_r)+1)\wt(\bar g_{r})} +
\sum_{r>i} \frac{\wt(h[r])}{\log^2(\wt(\bar h_r)+1)\wt(\bar
h_{r})} \right)F_2.
$$

By Lemma \ref{lemineq}, the latter sum does not exceed a constant
times $F_2$ that, in turn does not exceed a constant times
$\left\|\Sigma_2 \right\|^2 + \left\| \Sigma_3 \right\|^2 $.

Combining all these inequalities, we get

$$\left\| \pi(g)-\pi(h) \right\|^2\geq \left\| \Sigma_2 \right\|^2+
\left\| \Sigma_3 \right\|^2\geq \frac{\rho (\dist_G
(g,h))^{2}}{\log^4 (\dist_G(g,h)+1)}.$$

\medskip

Assume that $g,h$ are in case $\bfS$, so $f, f'\in M_k$ for some
$k\geq 1$. Theorem \ref{linelw} implies that $\diam (M_k,
\dist_k)\leq d\, \log v(k)$. Therefore by inequalities (\ref{nk})
$$
\rho \left( \dist_G\left( f, f' \right) \right)^2 \leq \rho \left(\mu_k\, \diam M_k \right)^2\leq \rho
(d \lambda_k \log v(k))^2\leq (d+1)^2 (\lambda_k)^2\ll (\mu_k)^2\, .
$$
The last but one inequality above follows by monotonicity and
sub-additivity of $\rho$, as well as by the equality $\lambda_k=
\rho (\lambda_k \log v(k))$.

Therefore
\begin{equation}\label{lowers2}\begin{array}{l}
\left\| \Sigma_2 \right\|^2 \gg \frac{\wt(\bar
g_{i+1})}{\log^2(\wt(\bar g_{i+1})+1)\wt(f)} \rho \left(
\dist_G\left(f, f' \right) \right)^2 + \frac{\wt(\bar
h_{i+1})}{\log^2(\wt(\bar h_{i+1})+1)\wt(f')} \rho \left(
\dist_G\left(f, f' \right) \right)^2\, .\end{array}
\end{equation}

An argument similar to the one in case $\bfB$ allows to obtain the
lower bound in this case also.
\endproof

Lemma \ref{lowerbound} completes the proof of Part \bfII of
Theorem \ref{main}.

Now let us prove Part \bfIII.

\begin{lemma} Consider the group $F$ (the free product of
finite groups $M_k$, $k\ge 1$) with the metric induced by the word
metric on $G$. Then $F$ has asymptotic dimension {one}.
\end{lemma}

\proof {It is a straightforward consequence of \cite[Theorem
2.1]{DS}, stating that the asymptotic dimension of a countable group
endowed with a proper left invariant metric is the supremum over the
asymptotic dimensions of its finitely generated subgroups. {Any
finitely  generated subgroup of $F$ is inside a free product of
finitely many $M_k$ (hence of asymptotic dimension $\le 1$},
according to \cite{BD1})}.
\endproof

Note that the asymptotic dimension of each $H_i$ is $1$ by
\cite{BD1}. It remains to use \cite[Corollary 24]{BD} and conclude
that the asymptotic dimension of $G$ is {at most $2$.}

Theorem \ref{main} is proved.

We conclude with some open questions.

\begin{qn}
{Does every finitely generated group have a Hilbert space
compression gap of the form $(\frac{f}{\log x}, f)$ for some
function $f : \R_+ \to \R_+$ ?}
\end{qn}

\begin{qn} Is there an amenable group with Hilbert space compression
$<\frac12$?
\end{qn}

\begin{qn} What is a \{uniformly convex Banach space\}--compression of  R.~Thompson's group $F$ or of the wreath product $\Z\wr\Z$?
\end{qn}





\end{document}